\documentclass[10pt]{article}
\usepackage{amsmath,amsfonts}
\usepackage{verbatim}
\usepackage{relsize}
\usepackage[inline]{enumitem}
\normalbaselines
\newtheorem{teor}{Theorem}

\newtheorem{lema}[teor]{Lemma}
\newtheorem{coro}[teor]{Corollary}
\newtheorem{rem}[teor]{Remark}

\newenvironment{demo}{\rm \trivlist \item[\hskip \labelsep{\it
      Proof}.]}{\nopagebreak \hfill $\square$ \endtrivlist}

\title{Constant mean curvature spacelike hypersurfaces in spatially open GRW spacetimes}

\author{Jos\'e A. S. Pelegr\'in${}^{a}$ and Marco Rigoli${}^{b}$ \\[6mm]
${}^a$ Departamento de Geometr\'\i a y Topolog\'\i a, \\ [0.5mm]
Universidad de Granada, 18071 Granada, Spain \\ E-mail\textup{:\texttt{jpelegrin@ugr.es}} \\[3mm]
${}^b$ Dipartimento di Matematica, \\[0.5mm] Universit\`a degli Studi di Milano, 20133 Milano, Italy,\\[0.5mm] E-mail\textup{: \texttt{marco.rigoli55@gmail.com}}\\[3mm]}

\date{}

\begin{document}

\maketitle

\thispagestyle{empty}

\begin{abstract}
In this paper we provide under certain geometric and physical
assumptions new uniqueness and non-existence
results for complete spacelike hypersurfaces of constant mean curvature in spatially open Generalized Robertson-Walker spacetimes. Some of our results are then applied to relevant spacetimes as the steady state spacetime, Einstein-de Sitter spacetime and certain radiation models.
\end{abstract}
\vspace*{5mm}



\section{Introduction}

Spacelike hypersurfaces of constant mean curvature in a spacetime are geometric objects of relevant physical and mathematical interest. They are critical points of the area functional under a suitable volume constraint \cite{brasil} and for instance, they play an important role in General Relativity as convenient initial data for the Cauchy problem \cite{CG, Li} and in the proof of the positivity of the gravitational mass \cite{SY}. A summary of the many reasons justifying their physical interest can be found in \cite{M-T}. As a consequence, the literature on the subject provides a number of articles dealing with existence and uniqueness results (see for instance \cite{ARS}, \cite{Ch} and \cite{PRR1}) in a large variety of spacetimes.

From a mathematical perspective, they exhibit nice Bernstein-type properties. In fact, the only complete spacelike hypersurfaces in the $(m+1)$-dimensional Lorentz-Minkowski spacetime $\mathbb{L}^{m+1}$ with zero mean curvature are spacelike hyperplanes. This striking result was proved by Calabi in \cite{Ca} for $m \leq 4$ and later extended to arbitrary dimension by Cheng and Yau in \cite{CY}. Their approach is based on a Simons-type formula for spacelike hypersurfaces in $\mathbb{L}^{m+1}$ and the application of a generalized maximum principle due to Omori \cite{Om} and Yau \cite{Y}.

In this article, we study constant mean curvature spacelike hypersurfaces in the family of cosmological models known as Generalized Robertson-Walker (GRW) spacetimes. They are warped products whose negative definite base represents a universal time and whose fiber is an arbitrary Riemannian manifold (see Section \ref{s2} for definitions). These spacetimes are relevant in physics and extend the classical notion of Robertson-Walker spacetime to the case when the fiber does not necessarily have constant sectional curvature \cite{ARS}. Thus, GRW spacetimes are not necessarily spatially homogeneous, making them suitable cosmological models to describe the universe in a more accurate scale \cite{RaS}. Furthermore, deformations of the metric on the fiber of classical Robertson-Walker spacetimes as well as GRW spacetimes with a time-dependent conformal change of metric also fit into the class of GRW spacetimes.

Classically, spacelike hypersurfaces have been studied in the class of GRW spacetimes known as spatially closed GRW spacetimes (see, for instance \cite{ARS} and \cite{PRR3}). However, some observations and theoretical arguments about the total mass balance of the universe \cite{Chiu} suggest the convenience of using spatially open models to describe our current universe. Moreover, spatially closed models also lead to a violation of the holographic principle \cite{B}. Due to these reasons, many authors have recently studied maximal hypersurfaces (i.e., spacelike hypersurface with zero mean curvature) in spatially open GRW spacetimes under certain other assumptions on the fiber, such as parabolicity \cite{RRS} or flatness \cite{PRR2}, and sometimes imposing some requirements, like the property of having a bounded hyperbolic angle or that of being contained in a slab \cite{RRS2}. 

The aim of this paper is to obtain new results for spacelike hypersurfaces of constant mean curvature in GRW spacetimes without any compactness assumption on the fiber and neither imposing a priori bounds on the hyperbolic angle nor the property of being contained in a slab. This will allow us to deal, respectively, with spacelike hypersurfaces that approach the null boundary as well as the future and past infinities. In order to obtain our results, we will assume a weak energy condition on the spacetime, the Null Energy Condition (NCC). Note that this energy condition is weaker than the Timelike Convergence Condition (see \cite{O'N}), which is a physically reasonable assumption commonly used to obtain uniqueness results for these hypersurfaces (see, for instance, \cite{BF} and \cite{N}).

Furthermore, as we point out in Remark \ref{remex}, our non-existence results have a clear physical meaning. Thus, we prove that in these ambient spacetimes there are no constant mean curvature spacelike hypersurfaces where the normal observers measure a different behaviour in the expansion/contraction of the universe than the comoving ones at every point.

\section{Preliminaries}
\label{s2} \noindent
Let $(F,g_{F})$ be an $m(\geq 2)$-dimensional (connected)
Riemannian manifold, $I$ an open interval in $\mathbb{R}$ endowed with the metric $-dt^2$  and $\rho$ a
positive smooth function defined on $I$. Then, the product manifold
$I \times F$ endowed with the Lorentzian metric
\begin{equation}\label{metrica}
\overline{g} = -\pi^*_{I} (dt^2) +\rho(\pi_{I})^2 \, \pi_{F}^* (g_{F})
\, ,
\end{equation}
where $\pi_{I}$ and $\pi_{F}$ denote the projections onto $I$ and
$F$, respectively, is called a GRW spacetime with
fiber $(F,g_{F})$, base $(I,-dt^2)$ and warping
function $\rho$. If the fiber has constant sectional curvature, it is called a Robertson-Walker spacetime.

In any GRW spacetime $\overline{M}=I\times_\rho F$, the
coordinate vector field $\partial_t:=\partial/\partial t$ is
(unitary) timelike, and hence $\overline{M}$ is time-orientable. On the other hand, if we consider the timelike vector field $K: =~ \rho({\pi}_I)\,\partial_t$, from the relation between the
Levi-Civita connection of $\overline{M}$ and those of the base and
the fiber \cite[Cor. 7.35]{O'N}, it follows that

\begin{equation}
\label{conexion} 
\overline{\nabla}_XK = \rho'({\pi}_I) X, 
\end{equation}

\noindent for any $X\in \mathfrak{X}(\overline{M})$, where $\overline{\nabla}$
is the Levi-Civita connection of the Lorentzian metric
(\ref{metrica}). Thus, $K$ is a closed conformal vector field.

From (\ref{conexion}) we easily see that the divergence on $\overline{M}$ of the reference frame $\partial_t$ satisfies 
$\mathrm{div}(\partial_t) = m \frac{\rho'(t)}{\rho(t)}$. Therefore, the observers in $\partial_t$ are spreading out (resp. coming 
together) if $\rho'>0$ (resp. $\rho'<0$).

Given an $m$-dimensional manifold $M$, an immersion $\psi: M
\rightarrow \overline{M}$ is said to be spacelike if the
Lorentzian metric (\ref{metrica}) induces, via $\psi$, a Riemannian
metric $g$ on $M$. In this case, $M$ is called a spacelike
hypersurface. We will denote by $\tau:=\pi_I\circ \psi$ the
restriction of $\pi_I$ along $\psi$.

The time-orientation of $\overline{M}$ allows to take, for each
spacelike hypersurface $M$ in $\overline{M}$, a unique unitary
timelike vector field $N \in \mathfrak{X}^\perp(M)$ globally defined
on $M$ with the same time-orientation of $\partial_t$, i.e., such
that $\overline{g}(N,\partial_t)\leq -1$. Note that $\overline{g}(N,\partial_t)=
-1$ at a point $p\in M$ if and only if $N = \partial_t$ at
$p$. We will denote by $A$ the shape operator associated to $N$.
Then, the mean curvature function associated to $N$ is given
by $H:= -(1/m) \mathrm{trace}(A)$. As it is well-known, the mean
curvature is constant if and only if the spacelike hypersurface is,
locally, a critical point of the $m$-dimensional area functional for
compactly supported normal variations, under certain constraints of
the volume. When the mean curvature vanishes identically, the
spacelike hypersurface is called a maximal hypersurface.

For a spacelike hypersurface $\psi: M \rightarrow \overline{M}$ with
Gauss map $N$, the hyperbolic angle $\varphi$, at any point
of $M$, between the unit timelike vectors $N$ and $\partial_t$, is
given by $\cosh \varphi=-\overline{g}(N,\partial_t)$. For simplicity,
throughout this paper we will refer to $\varphi$  as the hyperbolic angle function on $M$.

In any GRW spacetime $\overline{M}= I \times_\rho F$ there
is a remarkable family of spacelike hypersurfaces, namely its
spacelike slices $\{t_{0}\}\times F$, $t_{0}\in I$. It can be easily
seen that a spacelike hypersurface in $\overline{M}$ is a (piece of)
spacelike slice if and only if the function $\tau$ is constant.
Furthermore, a spacelike hypersurface in $\overline{M}$ is a (piece
of) spacelike slice if and only if the hyperbolic angle $\varphi$
vanishes identically. The shape operator of the spacelike slice
$\tau=t_{0}$ is given by $A=-\rho'(t_{0})/\rho(t_{0})\mathbb{I}$, where $\mathbb{I}$
denotes the identity transformation, and therefore its (constant)
mean curvature $H$ is given by $\rho'(t_{0})/\rho(t_{0})$. Thus, a spacelike slice
is maximal if and only if $\rho'(t_{0})=0$ (and hence, totally
geodesic).

\section{The set up} 
\label{sesu}
Let $\psi: M \rightarrow \overline{M}$ be an $m$-dimensional
spacelike hypersurface immersed in a GRW spacetime
$\overline{M}= I \times_\rho F$. Denoting by

$$\partial_t^T:= \partial_t+\overline{g}(N,\partial_t)N $$

\noindent the tangential component of $\partial_t$ along $\psi$, then it is
easy to check that the gradient of $\tau$ on $M$ is

\begin{equation}
\label{part}
\nabla \tau=-\partial_t^T
\end{equation}

\noindent and so

\begin{equation}
\label{sinh}
|\nabla \tau|^2=g(\nabla \tau,\nabla \tau)=\sinh^2 \varphi.
\end{equation}

Moreover, since the tangential
component of $K$ along $\psi$ is given by $K^T=K+\overline{g}(K,N)N$, a direct computation from
(\ref{conexion}) gives

\begin{equation}
\label{gradcosh}
\nabla \overline{g}(K,N)=-AK^T
\end{equation}

\noindent where we have used (\ref{part}), and also the fact that

\begin{equation}
\label{gch}
\nabla \cosh \varphi=A\partial_t^T + \frac{\rho'(\tau)}{\rho(\tau)} \cosh \varphi \ \partial_t^T.
\end{equation}

On the other hand, if we represent by $\nabla$ the Levi-Civita
connection of the metric $g$, then the Gauss and Weingarten
formulas for the immersion $\psi$ are respectively given by

\begin{equation}
\label{GF}
\overline{\nabla}_X Y=\nabla_X Y-g(AX,Y)N
\end{equation}

\noindent and

\begin{equation}
\label{WF}
AX=-\overline{\nabla}_X N,
\end{equation}

\noindent where $X,Y\in\mathfrak{X}({M})$. Then,  taking the tangential component in
(\ref{conexion}) and using (\ref{GF}) and (\ref{WF}), we obtain

\begin{equation}
\label{KT}
\nabla_X K^T=-\rho(\tau) \ \overline{g}(N,\partial_t) \ A X + \rho'(\tau) X,
\end{equation}

\noindent where $X\in\mathfrak{X}({M})$ and $\rho'(\tau):=\rho'\circ \tau$. Next, from (\ref{part}) we get

\begin{equation}
\label{nt}
\nabla_X \partial_t^T = \frac{\rho'(\tau)}{\rho(\tau)} \ g(X, \partial_t^T) \ \partial_t^T + \cosh \varphi \ A X + \frac{\rho'(\tau)}{\rho(\tau)} X.
\end{equation}

We now choose a local orthonormal reference frame $\{E_1, \dots, E_m \}$ on $TM$ to calculate the Laplacian of $\cosh \varphi$ as follows

\begin{equation}
\label{lap1}
\Delta \cosh \varphi  = \sum_{i=1}^m g(\nabla_{E_i} (A \partial_t^T), E_i) +   \sum_{i=1}^m g \left(\nabla_{E_i} \left(\frac{\rho'(\tau)}{\rho(\tau)} \cosh \varphi \ \partial_t^T \right), E_i\right),
\end{equation}

\noindent where we have made use of (\ref{gch}). After several computations we rewrite (\ref{lap1}) in the form

\begin{eqnarray}
\label{lap2}
\Delta \cosh \varphi &=& \sum_{i=1}^m g(\nabla_{E_i} (A \partial_t^T), E_i) -
 \frac{\rho''(\tau)}{\rho(\tau)} \cosh \varphi \sum_{i=1}^m g(E_i, \partial_t^T) g(\partial_t^T, E_i)  \\ \nonumber
& & + \frac{\rho'(\tau)^2}{\rho(\tau)^2} \cosh \varphi \sum_{i=1}^m g(E_i, \partial_t^T) g(\partial_t^T, E_i) + \frac{\rho'(\tau)}{\rho(\tau)} \sum_{i=1}^m g(E_i, A \partial_t^T) g(\partial_t^T, E_i) \\ \nonumber
& & + \frac{\rho'(\tau)^2}{\rho(\tau)^2} \cosh \varphi \sum_{i=1}^m g(E_i, \partial_t^T) g(\partial_t^T, E_i) + \frac{\rho'(\tau)}{\rho(\tau)} \cosh \varphi \sum_{i=1}^m g(\nabla_{E_i} \partial_t^T, E_i) \nonumber.
\end{eqnarray}

Taking into account (\ref{sinh}) and $(\nabla_X A)Y = \nabla_X (AY) - A(\nabla_X Y)$ for all $X, Y \in \mathfrak{X}(M)$, from (\ref{lap2}) we have

\begin{eqnarray}
\label{lap3}
\Delta \cosh \varphi &=& \sum_{i=1}^m g((\nabla_{E_i} A) \partial_t^T, E_i) + \sum_{i=1}^m g(\nabla_{E_i} \partial_t^T, A E_i) - \frac{\rho''(\tau)}{\rho(\tau)} \cosh \varphi \sinh^2 \varphi \\ \nonumber
& &  + 2 \frac{\rho'(\tau)^2}{\rho(\tau)^2} \cosh \varphi \sinh^2 \varphi + \frac{\rho'(\tau)}{\rho(\tau)} g(A \partial_t^T, \partial_t^T) + \frac{\rho'(\tau)}{\rho(\tau)} \cosh \varphi \sum_{i=1}^m g(\nabla_{E_i} \partial_t^T, E_i) \nonumber.
\end{eqnarray}

Furthermore, Codazzi equation $\overline{g}(\overline{\mathrm{R}}(X,Y) N, Z) = \overline{g}((\nabla_Y A_N) X, Z) - \overline{g}((\nabla_X A_N) Y, Z)$ yields

\begin{eqnarray}
\label{lap4}
\Delta \cosh \varphi &=& \sum_{i=1}^m \overline{g}(\overline{\mathrm{R}}(\partial_t^T, E_i)N, E_i) + \sum_{i=1}^m g((\nabla_{\partial_t^T} A) E_i, E_i) + \sum_{i=1}^m g(\nabla_{E_i} \partial_t^T, A E_i) \\ \nonumber 
& & - \frac{\rho''(\tau)}{\rho(\tau)} \cosh \varphi \sinh^2 \varphi  + 2 \frac{\rho'(\tau)^2}{\rho(\tau)^2} \cosh \varphi \sinh^2 \varphi + \frac{\rho'(\tau)}{\rho(\tau)} g(A \partial_t^T, \partial_t^T) \\ \nonumber
& & + \frac{\rho'(\tau)}{\rho(\tau)} \cosh \varphi \sum_{i=1}^m g(\nabla_{E_i} \partial_t^T, E_i) \nonumber.
\end{eqnarray}

Using (\ref{nt}) into (\ref{lap4}) we then obtain

\begin{eqnarray}
\label{lap5}
\Delta \cosh \varphi &=& \sum_{i=1}^m \overline{g}(\overline{\mathrm{R}}(\partial_t^T, E_i)N, E_i) + \sum_{i=1}^m g((\nabla_{\partial_t^T} A) E_i, E_i) + \frac{\rho'(\tau)}{\rho(\tau)} g(A \partial_t^T, \partial_t^T) \\ \nonumber 
& & + \cosh \varphi \ \mathrm{trace}(A^2) - m \frac{\rho'(\tau)}{\rho(\tau)} H  - \frac{\rho''(\tau)}{\rho(\tau)} \cosh \varphi \sinh^2 \varphi \\ \nonumber
& & + 2 \frac{\rho'(\tau)^2}{\rho(\tau)^2} \cosh \varphi \sinh^2 \varphi + \frac{\rho'(\tau)}{\rho(\tau)} g(A \partial_t^T, \partial_t^T)  + \frac{\rho'(\tau)^2}{\rho(\tau)^2} \cosh \varphi \sinh^2 \varphi \\ \nonumber
& & - m \frac{\rho'(\tau)}{\rho(\tau)} H \cosh^2 \varphi + m \frac{\rho'(\tau)^2}{\rho(\tau)^2} \cosh \varphi.
\end{eqnarray}

Since covariant derivations commute with contractions, choosing our local base in $T_pM$ satisfying $\left(\nabla_{E_j} E_i \right)_p = 0$ we obtain

 \begin{eqnarray}
\label{lap5}
\Delta \cosh \varphi &=& - \overline{\mathrm{Ric}}(\partial_t^T, N) - m \ g(\nabla H, \partial_t^T) + 2 \frac{\rho'(\tau)}{\rho(\tau)} g(A \partial_t^T, \partial_t^T) \\ \nonumber 
& & + \cosh \varphi \ \mathrm{trace}(A^2) - m \frac{\rho'(\tau)}{\rho(\tau)} H (\cosh^2 \varphi + 1)  - \frac{\rho''(\tau)}{\rho(\tau)} \cosh \varphi \sinh^2 \varphi \\ \nonumber
& & + 3 \frac{\rho'(\tau)^2}{\rho(\tau)^2} \cosh \varphi \sinh^2 \varphi + m \frac{\rho'(\tau)^2}{\rho(\tau)^2} \cosh \varphi.
\end{eqnarray}

Decomposing $N$ as $N=N_F-\overline{g}(N,\partial_t)\partial_t$, where
$N_F$ denotes the projection of $N$ on the fiber $F$, we know from \cite[Cor. 7.43]{O'N} that

\begin{equation}
\label{rict}
\overline{\mathrm{Ric}}(\partial_t,\partial_t)=-m \frac{\rho''(\tau)}{\rho(\tau)}
\end{equation}

and 

\begin{equation}
\label{ricNF}
\overline{\mathrm{Ric}}(N_F,N_F)= \mathrm{Ric}^F( N_F,N_F) + \sinh^2\varphi \left(\frac{\rho''(\tau)}{\rho(\tau)}+(m-1)\frac{\rho'(\tau)^2}{\rho(\tau)^2} \right),
\end{equation}

\noindent where $\mathrm{Ric}^F$ is the Ricci tensor of the fiber $F$. Therefore, from (\ref{rict}) and (\ref{ricNF}) we deduce

\begin{equation}
\label{ritn}
\overline{\mathrm{Ric}}(\partial_t^T, N) = -\cosh \varphi \left\{ \mathrm{Ric}^F( N_F,N_F) + (m-1) (\log \rho)''(\tau) \sinh^2 \varphi \right\}.
\end{equation}

We now insert (\ref{ritn}) into (\ref{lap5}) to obtain

\begin{eqnarray}
\label{lap6}
\Delta \cosh \varphi &=& \cosh \varphi \left\{ \mathrm{Ric}^F( N_F,N_F) - (m-1) (\log \rho)''(\tau) \sinh^2 \varphi \right\} - m \ g(\nabla H, \partial_t^T) \\ \nonumber 
& & + 2 \ \frac{\rho'(\tau)}{\rho(\tau)} g(A \partial_t^T, \partial_t^T) + \cosh \varphi \ \mathrm{trace}(A^2) - m \ \frac{\rho'(\tau)}{\rho(\tau)} H (\cosh^2 \varphi + 1)\\ \nonumber 
& &  - \frac{\rho''(\tau)}{\rho(\tau)} \cosh \varphi \sinh^2 \varphi + 3 \ \frac{\rho'(\tau)^2}{\rho(\tau)^2} \cosh \varphi \sinh^2 \varphi + m \ \frac{\rho'(\tau)^2}{\rho(\tau)^2} \cosh \varphi.
\end{eqnarray}

In the next step we compute $|\mathrm{Hess}(\tau)|^2$ using (\ref{nt}); having

\begin{eqnarray}
\label{he1}
|\mathrm{Hess}(\tau)|^2 &=& \sum_{i=1}^m g(\nabla_{E_i} \partial_t^T, \nabla_{E_i} \partial_t^T) = \frac{\rho'(\tau)^2}{\rho(\tau)^2} \sinh^4 \varphi + \cosh^2 \varphi \ \mathrm{trace}(A^2) \\ \nonumber
& & + m \ \frac{\rho'(\tau)^2}{\rho(\tau)^2} + 2 \ \frac{\rho'(\tau)}{\rho(\tau)} \cosh \varphi \ g(A \partial_t^T, \partial_t^T)+ 2 \ \frac{\rho'(\tau)^2}{\rho(\tau)^2} \sinh^2 \varphi \\ \nonumber 
& & - 2 \ m \ \frac{\rho'(\tau)}{\rho(\tau)} H \cosh \varphi.
\end{eqnarray}

We observe that (\ref{he1}) can be written as 

\begin{eqnarray}
\label{he2}
|\mathrm{Hess}(\tau)|^2 &=& \frac{\rho'(\tau)^2}{\rho(\tau)^2} \left( m - 1 + \cosh^ 4 \varphi \right) + \cosh^2 \varphi \ \mathrm{trace}(A^2) \\ \nonumber
& & + 2 \ \frac{\rho'(\tau)}{\rho(\tau)} \cosh \varphi \ g(A \partial_t^T, \partial_t^T) - 2 \ m \ \frac{\rho'(\tau)}{\rho(\tau)} H \cosh \varphi.
\end{eqnarray}

Using (\ref{lap6}) and (\ref{he2}), we have

\begin{eqnarray}
\label{clap1}
\cosh \varphi \ \Delta \cosh \varphi &=& \cosh^2 \varphi \left\{ \mathrm{Ric}^F( N_F,N_F) - (m-1) (\log \rho)''(\tau) \sinh^2 \varphi \right\} \\ \nonumber 
& & - m \cosh \varphi \ g(\nabla H, \partial_t^T) + |\mathrm{Hess}(\tau)|^2 \\ \nonumber
& & - \frac{\rho'(\tau)^2}{\rho(\tau)^2} \left( m - 1 + \cosh^ 4 \varphi \right)
+ 2 \ m \ \frac{\rho'(\tau)}{\rho(\tau)} H \cosh \varphi \\ \nonumber
& & - m \ \frac{\rho'(\tau)}{\rho(\tau)} H \cosh \varphi (\cosh^2 \varphi + 1)  - \frac{\rho''(\tau)}{\rho(\tau)} \cosh^2 \varphi \sinh^2 \varphi \\ \nonumber 
& & + 3 \ \frac{\rho'(\tau)^2}{\rho(\tau)^2} \cosh^2 \varphi \sinh^2 \varphi + m \ \frac{\rho'(\tau)^2}{\rho(\tau)^2} \cosh^2 \varphi.
\end{eqnarray}

Since $|\mathrm{Hess}(\tau)|^2 \geq 0$, after some computations (\ref{clap1}) leads to the differential inequality

\begin{eqnarray}
\label{clap2}
\cosh \varphi \ \Delta \cosh \varphi &\geq & \cosh^2 \varphi \left\{ \mathrm{Ric}^F( N_F,N_F) - m (\log \rho)''(\tau) \sinh^2 \varphi \right\} \\ \nonumber 
& & - m \cosh \varphi \ g(\nabla H, \partial_t^T) -m \ \frac{\rho'(\tau)}{\rho(\tau)} H \cosh \varphi \sinh^2 \varphi \\ \nonumber
& & + \frac{\rho'(\tau)^2}{\rho(\tau)^2} \sinh^2 \varphi (m + \sinh^2 \varphi).
\end{eqnarray}

Let us assume now that the ambient spacetime $\overline{M}$ satisfies the Null Convergence Condition (NCC). It is well known that a spacetime obeys the NCC if and only if its Ricci tensor satisfies $\overline{\mathrm{Ric}}(z,z) \geq 0$ for all lightlike vectors $z$. In particular, in an $(m+1)$-dimensional GRW spacetime NCC is equivalent to

\begin{equation}
\label{ncc}
\mathrm{Ric}^F(X_F,X_F) - m \ (\rho \rho'' - \rho'^2) g_F(X_F, X_F) \geq 0 ,
\end{equation}

\noindent for all $X_F$ tangent to the fiber $F$. This energy condition is a mathematical way to express that gravity, on average, attracts \cite{O'N} and it is automatically satisfied by any spacetime that obeys the Einstein field equations with a physically reasonable stress-energy tensor. Therefore, since $ \frac{1}{2} \Delta \sinh^2 \varphi = \cosh \varphi \ \Delta \cosh \varphi + |\nabla \cosh \varphi|^2$, from (\ref{clap2}) and (\ref{ncc}) we obtain

\begin{lema}
\label{lem1}
Let $\psi: M \longrightarrow \overline{M}$ be a spacelike hypersurface of constant mean curvature immersed in a GRW spacetime $\overline{M}= I \times_\rho F$ that obeys the NCC. Then, the hyperbolic angle of $M$ satisfies

\begin{equation}
\label{laps}
\frac{1}{2} \Delta \sinh^2 \varphi \geq -m \ \frac{\rho'(\tau)}{\rho(\tau)} H \cosh \varphi \sinh^2 \varphi + \frac{\rho'(\tau)^2}{\rho(\tau)^2} \sinh^2 \varphi (m + \sinh^2 \varphi).
\end{equation}

\end{lema}

\section{Main results}
\label{semr}

As a first consequence of Lemma \ref{lem1} we obtain the following non-existence result

\begin{teor}
\label{teo1}
Let $\overline{M}= I \times_\rho F$ be a GRW spacetime obeying the NCC. Then there are no complete spacelike hypersurfaces $\psi: M \longrightarrow \overline{M}$ with constant mean curvature $H$ satisfying

$$H \rho'(\tau) \leq 0 \ \ \ \ \ \inf_M \frac{\rho'(\tau)^2}{\rho(\tau)^2} >0$$

\noindent and the volume growth condition

\begin{equation}
\label{vg}
\liminf\limits_{r\rightarrow +\infty} \frac{\log(\mathrm{Vol}(B_r))}{r^2} < + \infty.
\end{equation}
\end{teor}

\begin{demo}
By contradiction, let us suppose the existence of such a spacelike hypersurface. From Lemma \ref{lem1} and the assumption $H \rho'(\tau) \leq 0$ we deduce that the hyperbolic angle of $M$ satisfies

\begin{equation}
\label{lsh}
\frac{1}{2} \Delta \sinh^2 \varphi \geq \frac{\rho'(\tau)^2}{\rho(\tau)^2} \sinh^2 \varphi (m + \sinh^2 \varphi).
\end{equation}

Completeness and the assumption on the volume growth of geodesic balls enables us to apply \cite[Thm. 4.2]{AMR} to the above differential inequality to infer that $\sinh^2 \varphi$ is bounded from above. Next, thanks to the assumptions $\inf_M \frac{\rho'(\tau)^2}{\rho(\tau)^2} >0$ and (\ref{vg}) we can use \cite[Thm. 4.1]{AMR} to deduce that the hyperbolic angle identically vanishes on $M$. It follows that $\psi(M)$ is a spacelike slice $\{t_0 \} \times F$ with mean curvature $H =\frac{\rho'(t_{0})}{\rho(t_{0})}$. However, by assumption $H \rho'(\tau) \leq 0$, so we immediately deduce that $\rho'(t_0) = 0$. This contradicts the assumption $\inf_M \frac{\rho'(\tau)^2}{\rho(\tau)^2} >0$.
\end{demo}

\begin{rem}
\normalfont
The conclusions of Theorem \ref{teo1} still hold substituting (\ref{vg}) by the integral condition

\begin{equation}
\label{inc}
\liminf\limits_{r\rightarrow +\infty} \frac{1}{r^2} \int_{B_r} |\sinh \varphi|^{p} < + \infty
\end{equation}

\noindent for some $p >0$ (see \cite[Rem. 4.2]{AMR}).
\end{rem}

The above technique can also be used to give a bound for the hyperbolic angle in the next theorem

\begin{teor}
\label{teo2}
Let $\overline{M}= I \times_\rho F$ be a GRW spacetime obeying the NCC and let $\psi: M \longrightarrow \overline{M}$ be a complete spacelike hypersurface with constant mean curvature $H$ satisfying the volume growth condition (\ref{vg}). If $0 < H \leq \frac{\rho'(\tau)}{\rho(\tau)}$ (resp., $\frac{\rho'(\tau)}{\rho(\tau)} \leq H < 0$) on $M$, then the hyperbolic angle of the immersion $\varphi$ satisfies 

$$ \sinh^2 \varphi \leq m (m-2).$$
\end{teor}

\begin{demo}
By means of Lemma \ref{lem1} and the condition $0 < H \leq \frac{\rho'(\tau)}{\rho(\tau)}$  (resp., $\frac{\rho'(\tau)}{\rho(\tau)} \leq H < 0$) we obtain the differential inequality

\begin{equation}
\label{lsh2}
\frac{1}{2} \Delta \sinh^2 \varphi \geq \frac{\rho'(\tau)^2}{\rho(\tau)^2} \sinh^2 \varphi \left( m -m \sqrt{1 + \sinh^2 \varphi} + \sinh^2 \varphi \right).
\end{equation}

\noindent where $\inf_M \frac{\rho'(\tau)^2}{\rho(\tau)^2} >0$ because $H >0$. From (\ref{vg}) and \cite[Thm. 4.2]{AMR} we thus obtain that $\sinh^2 \varphi$ is bounded from above. Again by $\inf_M \frac{\rho'(\tau)^2}{\rho(\tau)^2} >0$ and (\ref{vg}), we can apply \cite[Thm. 4.1]{AMR} to obtain 

$$\sup_M(\sinh^2 \varphi) \left( m -m \sqrt{1 + \sup_M(\sinh^2 \varphi)} + \sup_M(\sinh^2 \varphi) \right) \leq 0,$$

\noindent that gives the desired conclusion.
\end{demo}

In particular, for spacelike surfaces immersed in a $3$-dimensional GRW spacetime, we have

\begin{coro}
\label{cordim2}
Let $\overline{M}= I \times_\rho F$ be a GRW spacetime of dimension $3$ obeying the NCC and let $\psi: M \longrightarrow \overline{M}$ be a complete spacelike hypersurface with constant mean curvature $H$ satisfying (\ref{vg}). If $0 < H \leq \frac{\rho'(\tau)}{\rho(\tau)}$ (resp., $\frac{\rho'(\tau)}{\rho(\tau)} \leq H < 0$), then $M$ is a spacelike slice.
\end{coro}

For the next result we shall use a different technique. Note that in Remark \ref{reri} after the proof we give a sufficient condition for the existence of a function $\zeta$ satisfying $\textit{(i)}$, $\textit{(ii)}$ and $\textit{(iii)}$ below. 

\begin{teor}
\label{teodiv}
Let $\overline{M}= I \times_\rho F$ be a GRW spacetime obeying the NCC and let $\psi: M \longrightarrow \overline{M}$ be a complete spacelike hypersurface with constant mean curvature. Assume the existence of $\zeta \in C^0 (M\setminus B_{R_0})\cap Lip_{loc} (M \setminus B_{R_0})$ satisfying

\begin{center}
\begin{enumerate*}[label=(\roman*),itemjoin={,\quad}]
\item $\Delta \zeta \leq \Lambda_0 + \Lambda_1 \zeta$
\item $|\nabla \zeta| \leq \Lambda_2$
\item $\zeta(x) \longrightarrow +\infty$ as $x \longrightarrow \infty$
\end{enumerate*}
\end{center}

\noindent on $M \setminus B_{R_0}$ for some positive constants $\Lambda_0, \Lambda_1, \Lambda_2$. Furthermore, suppose that $\rho'(\tau) H \leq 0$ on $M$ and that $\sinh \varphi \in~L^2(M)$. Then,

$$\rho'(\tau) \sinh^2 \varphi \equiv 0 \ \text{on} \ M.$$

\end{teor}

\begin{demo}
Fix $R > T > R_0 + 2$ and let $\alpha \in C^1 (\mathbb{R}_0^+) \cap C^2([0, R])$ be such that

$$ \alpha(t) \geq 0 \ \text{on} \ \mathbb{R}_0^+, \ \alpha(t)=1 \ \text{on} \ [0, T], \ \alpha(t)=0 \ \text{on} \ [R, + \infty)$$

and

\begin{equation}
\label{alf}
\alpha'(t) \leq \frac{P}{R-T}, \ \ \ |\alpha''(t)| \leq \frac{P}{(R-T)^2}
\end{equation}

\noindent on $[0,R]$ for some constant $P>0$ independent of $R$ and $T$. Next, we define the cut-off function $\phi$ by setting

$$\phi = \begin{cases} \alpha(\zeta(x)) & \text{on} \ M\setminus D_{R_0} \\
\\
1 & \text{on} \ D_{R_0} \end{cases}$$

\noindent where $D_T = \{ x\in M : \zeta(x) < T\}$. We then have 

\begin{equation}
\label{psi}
\begin{cases} \phi(x) \geq 0 \ \text{on} \ M, \ \phi(x)=1 \ \text{on} \ D_T, \ \phi(x)=0 \ \text{on} \ M \setminus D_R \\
\\
\nabla \phi = 0 \ \text{on} \ \partial D_R \cup D_T \cup (M \setminus D_R), \ \nabla \phi \leq \frac{P}{R-T} \end{cases}
\end{equation}

Finally, from our assumptions on $\zeta$ and (\ref{alf}) we deduce

\begin{equation}
\label{lapsi}
\Delta \phi = 0 \ \text{on} \ D_T \cup (M \setminus D_R) , \ \Delta \phi \leq C \frac{R}{R-T} + \frac{D}{R-T} + \frac{E}{(R-T)^2} \ \text{on} \ \overline{D}_R \setminus D_T
\end{equation}

\noindent for some positive constants $C, D, E$ independent of $R$ and $T$. Next we consider any $C^2$ function $u \geq 0$ on $M$. By the divergence theorem and the properties of $\phi$ we have

\begin{equation}
\label{i1}
0 = \int_{D_R} \phi \ \Delta u \ dV + \int_{D_R \setminus D_T} g(\nabla \phi, \nabla u) \ dV,
\end{equation} 

\noindent with $dV$ the volume element of $M$. Observing that

$$ \mathrm{div}(u \ \nabla \phi) = u \ \Delta \phi + g(\nabla u, \nabla \phi)$$

\noindent and that $\nabla \phi= 0$ on $\partial D_R \cup \partial D_T$, a further application of the divergence theorem yields

\begin{equation}
\label{i2}
0 = \int_{D_R \setminus D_T} u \ \Delta \phi \ dV + \int_{D_R \setminus D_T} g(\nabla u, \nabla \phi) \ dV.
\end{equation} 

\noindent Inserting (\ref{i2}) into (\ref{i1}) we deduce

\begin{equation}
\label{i3}
\int_{D_R} \phi \ \Delta u \ dV = \int_{D_R \setminus D_T} u \ \Delta \phi \ dV.
\end{equation} 

We set $T = \frac{R}{2}$ and using (\ref{lapsi}) we finally obtain

\begin{equation}
\label{i4}
\int_{D_{R/2}} \Delta u \ dV \leq G \int_{D_R \setminus D_{R/2}} u \ dV
\end{equation}

\noindent for some constant $G>0$ sufficiently large and independent of $R$. Now, let $u=\sinh^2 \varphi$ and use $\rho'(\tau) H \leq 0$ to deduce from (\ref{laps}) of Lemma \ref{lem1}

\begin{equation}
\label{inequ}
\frac{1}{2} \Delta u \geq \frac{\rho'(\tau)^2}{\rho(\tau)^2} u \ (m+u).
\end{equation}

\noindent Inserting (\ref{inequ}) into (\ref{i4}), we get

\begin{equation}
\label{i5}
\int_{D_{R/2}} \frac{\rho'(\tau)^2}{\rho(\tau)^2} u \ (m+u) \ dV \leq \frac{G}{2} \int_{D_R \setminus D_{R/2}} u \ dV.
\end{equation}

\noindent Since $u \in L^1(M)$ and $\zeta(x) \longrightarrow +\infty$ as $x \longrightarrow \infty$, letting $R \longrightarrow +\infty$ in (\ref{i5}) we conclude that

$$\rho'(\tau) \sinh^2 \varphi \equiv 0 \ \text{on} \ M.$$
\end{demo}

\begin{rem}
\label{reri}
\normalfont
Let $r(x)= \mathrm{dist}(x, o)$ for some fixed origin $o$ in the complete manifold $M$. Observe that, by the Laplacian comparison theorem, the condition

\begin{equation}
\label{ricr}
\mathrm{Ric}(\nabla r, \nabla r) \geq -(m-1)B^2 (1+r^2) \ \text{on} \ M
\end{equation}

\noindent for some constant $B>0$ implies

\begin{equation}
\label{dr}
\Delta r \leq \Lambda_0 + \Lambda_1 r \ \text{on} \ M \setminus B_1
\end{equation} 

\noindent for some positive constants $\Lambda_0, \Lambda_1$ and obviously by Gauss lemma $|\nabla r| = 1$ (see for instance \cite{BMR}). This observation yields the following
\end{rem}

\begin{teor}
\label{teorib}
Let $\overline{M}= I \times_\rho F$ be a GRW spacetime with fiber $F$ of non-negative sectional curvature and warping function $\rho$ satisfying $(\log \rho)'' \leq 0$.  Let $\psi: M \longrightarrow \overline{M}$ be a complete spacelike hypersurface of constant mean curvature $H$ and suppose that $\rho'(\tau) H \leq 0$ on $M$ and $\sinh \varphi \in L^2(M)$. Then,
$$\rho'(\tau) \sinh^2 \varphi \equiv 0 \ \text{on} \ M.$$

\end{teor}

\begin{demo}
Using \cite[Lemma 13]{ARR} and the assumptions of the theorem we deduce that the Ricci curvature of $M$ is bounded from below. Moreover, from our assumptions the spacetime obeys the NCC. Therefore, we see from Remark \ref{reri} that we can use Theorem \ref{teodiv} choosing as $\zeta$ the distance function $r(x)$ on $M$ to prove the desired conclusion.  
\end{demo}

\begin{rem}
\normalfont
Note that the existence of $\zeta$ satisfying conditions \textit{(i)} and \textit{(iii)} in Theorem \ref{teodiv} implies the validity of the weak maximum principle for the Laplacian on the complete manifold $M$ \cite[Thm. 3.1]{AMR}. Observe that the conclusion of Theorem \ref{teodiv} can also be expressed in the form

\begin{equation}
\label{sush}
\mathrm{supp} \ \varphi \subseteq \{\tau \in I: \rho'(\tau) = 0 \}.
\end{equation}

\end{rem}

A different integrability request on $\sinh^2 \varphi$ enables us to obtain conclusion (\ref{sush}) avoiding the assumptions on the function $\zeta$ of Theorem \ref{teodiv}. Indeed, this is the content of the next result.

\begin{teor}
\label{teo3}
Let $\overline{M}= I \times_\rho F$ be a GRW spacetime obeying the NCC and let $\psi: M \longrightarrow \overline{M}$ be a complete spacelike hypersurface of constant mean curvature $H$. If $H \rho'(\tau) \leq 0$ and $\sinh^2 \varphi \in L^q (M)$ for some $q >2$, then

$$ \mathrm{supp} \ \varphi \subseteq \{\tau \in I: \rho'(\tau) = 0 \}.$$
\end{teor}

\begin{demo}
From Lemma \ref{lem1}, (\ref{laps}) holds on $M$ and, setting $u = \sinh^2 \varphi$, we have

\begin{equation}
\label{lu1}
u \Delta u - 2 \ m \ \frac{\rho'(\tau)^2}{\rho(\tau)^2} u^2 - 2 \ \frac{\rho'(\tau)^2}{\rho(\tau)^2} u^3 \geq 0.
\end{equation}

\noindent Consider the operator

$$L = \Delta -  m \ q \ \frac{\rho'(\tau)^2}{\rho(\tau)^2}$$

\noindent and observe that any positive constant $\kappa$ solves $L \kappa \leq 0$ on $M$. Using now \cite[Thm. 3.3]{MRS} with the choices (in the notation of Theorem 3.3) $H=\frac{1}{2} q$, $\beta = \frac{1}{2} q -1$, $K=0$, $p=2$ and  taking into account that the assumption $u \in L^q (M)$ for some $q >2$ implies 

$$\frac{1}{\int_{\partial B_r} u^{2(\beta+1)}} \notin L^1 (+\infty),$$

\noindent (see \cite{RS}) we get that (\ref{lu1}) has no non-negative $C^2$-solutions $u$ on $M$ satisfying

$$ \mathrm{supp} \ u \cap \left\{\tau \in I: 2 \ \frac{\rho'(\tau)^2}{\rho(\tau)^2} > 0 \right\} = \emptyset.$$
\end{demo}

As a consequence of Theorem \ref{teo3} we have

\begin{coro}
\label{co}
Let $\overline{M}= I \times_\rho F$ be a GRW spacetime obeying the NCC. Then, there are no complete spacelike hypersurfaces $M$ of constant mean curvature in $\overline{M}$ such that $H \rho'(\tau) \leq 0$, $\sinh^2 \varphi \in L^q (M)$ for some $q >2$ and $\rho'(\tau) \neq 0$ on $M$.
\end{coro}

In what follows we prove a result that leads to interesting consequences for certain well-known spacetimes.

\begin{teor}
\label{teoale}
Let $\overline{M}= I \times_\rho F$ be a GRW spacetime whose warping function satisfies $(\log \rho)'' \leq 0$ and whose fiber $F$ has non-negative sectional curvature. Then, there are no complete spacelike hypersurfaces of constant mean curvature in $\overline{M}$ satisfying $H \rho'(\tau) \leq 0$ and $\inf_M \frac{\rho'(\tau)^2}{\rho(\tau)^2} >0$.
\end{teor}

\begin{demo}
Let us suppose the existence of such a spacelike hypersurface $\psi: M \longrightarrow \overline{M}$. From Lemma \ref{lem1} and our assumptions, that imply the NCC, the hyperbolic angle of the immersion verifies

\begin{equation}
\label{ni}
\frac{1}{2} \Delta \sinh^2 \varphi \geq \inf_M \left(\frac{\rho'(\tau)^2}{\rho(\tau)^2}\right) \sinh^4 \varphi.
\end{equation}

Moreover, from \cite[Lemma 13]{ARR} the Ricci curvature of $M$ is bounded from below. Since $M$ is complete and $\inf_M \frac{\rho'(\tau)^2}{\rho(\tau)^2} >0$, we can use \cite[Lemma 2]{PRR2} (this is a result obtained in \cite{N}) on (\ref{ni}) to conclude that the hyperbolic angle identically vanishes on $M$. Now, from the assumption $H \rho'(\tau) \leq 0$ we deduce that $M$ is a totally geodesic spacelike slice, contradicting $\inf_M \frac{\rho'(\tau)^2}{\rho(\tau)^2} >0$.
\end{demo}

\begin{rem}
\normalfont
Note that we can also guarantee the boundedness of the Ricci curvature of the spacelike hypersurface $M$ if the fiber has Ricci curvature bounded from below (not necessarily by zero) and the hyperbolic angle of $M$ is bounded. Indeed, the aim of this bound on the Ricci curvature of $M$ is to guarantee that the Omori-Yau maximum principle for the Laplacian holds on $M$. Even more, we can substitute this assumption by requiring a controlled decay of the Ricci curvature \cite{AMR}.
\end{rem}

\begin{rem}
\label{remex}
\normalfont
We give a physical interpretation to the assumptions in our theorems referring to \cite{SW2}. In order to do so, at each point $p \in M$ we define in a neighborhood $U$ of $p$ in the spacetime  a unitary future-pointing timelike vector field $\widetilde{N}$ such that $\widetilde{N} = N$ on $U \cap M$, being $N$ the unit normal vector field to $M$. If we compute the divergence in $\overline{M}$ of $\widetilde{N}$ at $p \in M$ we get

\begin{equation}
\label{diin}
\overline{\mathrm{div}}(\widetilde{N})_p = m H(p),
\end{equation}

\noindent where $H$ is the mean curvature function of $M$. Since the integral curves of this vector field $\widetilde{N}$ are known as the normal observers, if the mean curvature $H(p)$ is positive (resp. negative) at some point $p \in M$, these normal observers will measure that they are spreading out (resp. coming together).

Moreover, in a GRW spacetime $\overline{M}= I \times_\rho F$ there is a distinguished family of observers known as the comoving observers, which are defined as the integral curves of the vector field $\partial_t$. Since the divergence of this vector field in $\overline{M}$ is

\begin{equation}
\label{diin}
\overline{\mathrm{div}}(\partial_t) = m \frac{\rho'}{\rho},
\end{equation}

\noindent we obtain that the comoving observers will measure that the spacetime is expanding or contracting depending on the sign of $\rho'$. Thus, our results in Theorem \ref{teo1}, Corollary \ref{cordim2} and Theorem \ref{teoale} imply that, under certain assumptions, there are no complete spacelike hypersurfaces of constant mean curvature in these ambient spacetimes where the normal observers measure that the universe is non-expanding whereas the comoving ones measure that it is non-contracting at some point $p \in M$ or vice versa.
\end{rem}

From Theorem \ref{teoale} we obtain the following non-existence results, which extend \cite[Corollaries 5, 6, 7]{PRR2} to the case of constant mean curvature spacelike hypersurfaces.

\begin{coro}
\label{ste}
There are no complete spacelike hypersurfaces of non-positive constant mean curvature in the $(m+1)$-dimensional steady state spacetime $\mathbb{R}\times_{e^t} \mathbb{R}^m$.
\end{coro}

\begin{coro}
\label{eds}
There are no complete spacelike hypersurfaces of non-positive constant mean curvature bounded away from future infinity in the $(m+1)$-dimensional Einstein-de Sitter spacetime $\mathbb{R}^+ \times_{t^{2/3}} \mathbb{R}^m$.
\end{coro}

\begin{coro}
\label{rad}
There are no complete spacelike hypersurfaces of non-positive constant mean curvature bounded away from future infinity in the $(m+1)$-dimensional Robertson-Walker radiation model $\mathbb{R}^+ \times_{(2at)^{1/2}} \mathbb{R}^m$, where $a>0$.
\end{coro}

\begin{demo}
We recall that for a spacelike hypersurface in a GRW spacetime, being bounded away from future infinity analytically means that $\sup_M \tau < + \infty$.
\end{demo}

\section*{Acknowledgements} The first author is supported by Spanish MINECO and ERDF project MTM2016-78807-C2-1-P.

\end{document}